\documentclass[12pt,british]{article}

\usepackage{babel}
\usepackage{amsmath,amsfonts,amssymb,amsthm}
\usepackage[latin1]{inputenc}
\usepackage[T1]{fontenc}
\usepackage{enumerate}
\usepackage{textcomp}
\usepackage{eucal}

\usepackage{epsfig}
\usepackage[noBBpl]{mathpazo}

\usepackage[matrix,arrow]{xy}

\theoremstyle{plain}
\newtheorem{thm}{Theorem}

\newtheorem{prop}[thm]{Proposition}
\newtheorem{cor}[thm]{Corollary}

\theoremstyle{remark}

\newcommand{\reqref}[1]{(\protect\ref{eq:#1})}
\newcommand{\resec}[1]{Section~\protect\ref{sec:#1}}
\newcommand{\repr}[1]{Proposition~\protect\ref{prop:#1}}

\newcommand{\reth}[1]{Theorem~\protect\ref{th:#1}}
\newcommand{\req}[1]{equation~(\protect\ref{eq:#1})}

\makeatletter

\renewcommand{\p@enumii}{}

\def\@enum@{\list{\csname label\@enumctr\endcsname}%
           {\usecounter{\@enumctr}\def\makelabel##1{
\normalfont\ignorespaces\emph{{##1}~}}
\setlength{\labelsep}{3pt}
\setlength{\parsep}{0pt}
\setlength{\itemsep}{0pt}
\setlength{\leftmargin}{0pt}
\setlength{\labelwidth}{0pt}
\setlength{\listparindent}{\parindent}
\setlength{\itemsep}{0pt}
\setlength{\itemindent}{0pt}
\topsep=3pt plus 1pt minus 1 pt}}

\def\@map#1#2[#3]{\mbox{$#1 \colon\thinspace #2 \to #3$}}
\def\map#1#2{\@ifnextchar [{\@map{#1}{#2}}{\@map{#1}{#2}[#2]}}

\makeatother

\newcommand{\rp}{\ensuremath{\mathbb{R}P^2}}
\newcommand{\dt}{\ensuremath{\mathbb D}^{2}}
\newcommand{\Z}{\ensuremath{\mathbb Z}}
\newcommand{\N}{\ensuremath{\mathbb N}}
\newcommand{\St}[1][2]{\ensuremath{\mathbb S}^{#1}}
\newcommand{\Et}{\ensuremath{\mathbb E}^{\,2}}
\newcommand{\FF}{\ensuremath{\mathbb F}}
\newcommand{\F}[1][n]{\ensuremath{\FF_{{#1}}}}
\newcommand{\comm}{\ensuremath{\rightleftharpoons}}
\renewcommand{\to}{\ensuremath{\longrightarrow}}
\renewcommand{\ker}[1]{\ensuremath{\operatorname{\text{Ker}}\left({#1}\right)}}
\DeclareRobustCommand*{\up}[1]{\textsuperscript{#1}}
\newcommand{\brak}[1]{\ensuremath{\left\{ #1 \right\}}}

\newcommand{\lhra}{\mathrel{\lhook\joinrel\to}}

\newcommand{\ang}[1]{\ensuremath{\left\langle #1\right\rangle}}

\newcommand{\setangr}[2]{\ensuremath{\ang{#1 \,\left\lvert \, #2 \right.}}}

\newcommand{\setr}[2]{\ensuremath{\brak{#1 \,\left\lvert \, #2 \right.}}}

\renewcommand{\epsilon}{\varepsilon}
\renewcommand{\th}{\ensuremath{\up{th}}}
\newcommand{\pnm}[1][n]{\ensuremath{P_{{#1}}(M)}}
\newcommand{\sn}[1][n]{\ensuremath{S_{{#1}}}}

\newcommand{\indalp}[2]{\ensuremath{\alpha_{#1,#2}}}
\newcommand{\indbet}[3]{\ensuremath{\beta_{#1,#2,#3}}}

\begin{document}

\title{The braid groups of the projective plane and the Fadell-Neuwirth short exact sequence}  

\author{Daciberg~Lima~Gon\c{c}alves\\
Departamento de Matem\'atica - IME-USP,\\
Caixa Postal~\textup{66281}~-~Ag.~Cidade de S\~ao Paulo,\\ CEP:~\textup{05311-970} - S\~ao Paulo - SP - Brazil.\\ e-mail: \textup{\texttt{dlgoncal@ime.usp.br}}\vspace*{4mm}\\
John~Guaschi\\
Laboratoire de Math\'ematiques Emile Picard,\\ 
Institut de Math\'ematiques de Toulouse UMR CNRS~\textup{5219},\\ 
UFR-MIG, Universit\'e Toulouse~III,\\ \textup{31062}~Toulouse Cedex~\textup{9}, France.\\
e-mail: \textup{\texttt{guaschi@picard.ups-tlse.fr}}}

\date{12\up{th}~April~2007}           

\maketitle

\begin{abstract}
\noindent 
We study the pure braid groups $P_n(\rp)$ of the real projective plane
$\rp$, and in particular the possible splitting of the Fadell-Neuwirth
short exact sequence $1 \to P_m(\rp \setminus\brak{x_1,\ldots,x_n})
\lhra P_{n+m}(\rp) \stackrel{p_{\ast}}{\to} P_n(\rp) \to 1$, where
$n\geq 2$ and $m\geq 1$, and $p_{\ast}$ is the homomorphism which
corresponds geometrically to forgetting the last $m$ strings. This
problem is equivalent to that of the existence of a section for the
associated fibration $\map{p}{F_{n+m}(\rp)}[F_n(\rp)]$ of
configuration spaces. Van Buskirk proved in 1966 that $p$ and
$p_{\ast}$ admit a section if $n=2$ and $m=1$. Our main result in this
paper is to prove that there is no section if $n\geq 3$. As a
corollary, it follows that $n=2$ and $m=1$ are the only values for
which a section exists. As part of the proof, we derive a presentation
of $P_n(\rp)$: this appears to be the first time that such a presentation has been given in the literature.

\end{abstract}

\section{Introduction}\label{sec:intro}

Braid groups of the plane were defined by Artin in~1925~\cite{A1}, and
further studied in~\cite{A2,A3}. They were later generalised using the
following definition due to Fox~\cite{FoN}. Let $M$ be a compact,
connected surface, and let $n\in\N$. We denote the set of all ordered
$n$-tuples of distinct points of $M$, known as the \emph{$n\th$
configuration space of $M$}, by: 
\begin{equation*}
F_n(M)=\setr{(p_1,\ldots,p_n)}{\text{$p_i\in M$ and $p_i\neq p_j$ if
$i\neq j$}}.
\end{equation*}
Configuration spaces play an important r\^ole in several branches of
mathematics and have been extensively studied, see~\cite{CG,FH} for
example. 

The symmetric group $\sn$ on $n$ letters acts freely on $F_n(M)$ by
permuting coordinates. The corresponding quotient will be denoted by
$D_n(M)$. Notice that $F_n(M)$ is a regular covering of $D_n(M)$. The
\emph{$n\th$ pure braid group $P_n(M)$} (respectively the \emph{$n\th$
braid group $B_n(M)$}) is defined to be the fundamental group of
$F_n(M)$ (respectively of $D_n(M)$). If $m\in \N$, then we may define
a homomorphism \( \map{p_{\ast}}{\pnm[n+m]}[\pnm] \) induced by the
projection \( \map{p}{F_{n+m}(M)}[F_n(M)] \) defined by \(
p((x_1,\ldots, x_n, \ldots, x_{n+m}))= (x_1,\ldots,x_n) \).
Representing $\pnm[n+m]$ geometrically as a collection of $n+m$
strings, $p_{\ast}$ corresponds to forgetting the last $m$ strings.
\textbf{We adopt the convention, that unless explicitly stated, all
homomorphisms \( \pnm[n+m] \to \pnm \) in the text will be this one}.

If $M$ is without boundary, Fadell and Neuwirth study the map $p$, and
show (\cite[Theorem~3]{FaN}) that it is a locally-trivial fibration.
The fibre over a point $(x_1,\ldots,x_n)$ of the base space is
$F_m(M\setminus\brak{x_1,\ldots,x_n})$ which we consider to be a
subspace of the total space via the map $\map
i{F_m(M\setminus\brak{x_1,\ldots,x_n})}[F_n(M)]$ defined by
$i((y_1,\ldots,y_m))=(x_1,\ldots,x_n,y_1,\ldots,y_m)$. Applying the
associated long exact sequence in homotopy, we obtain the \emph{pure
braid group short exact sequence of Fadell and Neuwirth}:
\begin{equation}\label{eq:split}
1 \to P_m(M\setminus\brak{x_1,\ldots,x_n}) \stackrel{i_{\ast}}{\to} \pnm[n+m] \stackrel{p_{\ast}}{\to}
\pnm \to 1, \tag{\normalfont\textbf{PBS}}
\end{equation}
where $n\geq 3$ if $M$ is the sphere $\St$~\cite{Fa,FvB}, $n\geq 2$ if $M$ is
the real projective plane $\rp$~\cite{vB}, and $n\geq 1$ otherwise~\cite{FaN},
and where $i_{\ast}$ and $p_{\ast}$ are the homomorphisms induced by the maps
$i$ and $p$ respectively. The sequence also exists for the classical pure braid
group $P_n$, where $M$ is the $2$-disc $\dt$ (or the plane). The short exact
sequence~\reqref{split} has been widely studied, and may be employed for example
to determine presentations of $P_n(M)$ (see \resec{pres}), its centre, and
possible torsion.  It was also used in recent work on the structure of the
mapping class groups~\cite{PR} and on Vassiliev invariants for surface
braids~\cite{GMP}. 

The decomposition of $P_n$ as a repeated semi-direct product of free groups (known as the `combing' operation) is the principal result of Artin's classical theory of braid groups~\cite{A2}, and allows one to obtain normal forms and to solve the word problem. More recently, it was used by Falk and Randell to study the lower central series and the residual nilpotence of $P_n$~\cite{FR}, and by Rolfsen and Zhu to prove that $P_n$ is bi-orderable~\cite{RZ}.

The problem of deciding whether such a decomposition exists for
surface braid groups is thus fundamental. This was indeed a recurrent
and central question during the foundation of the theory and its
subsequent development during the 1960's~\cite{Fa,FaN,FvB,vB,Bi1}. If
the fibre of the fibration is an Eilenberg-MacLane space then the
existence of a section for $p_{\ast}$ is equivalent to that of a
cross-section for $p$~\cite{Ba,Wh1} (cf.~\cite{GG3}). But with the
exception of the construction of sections in certain cases (for the
sphere~\cite{Fa} and the torus~\cite{Bi1}), no progress on the
possible splitting of~\reqref{split} was recorded for nearly forty
years. In the case of orientable surfaces without boundary of genus at
least two, the question of the splitting of~\reqref{split} which was
posed explicitly by Birman in 1969~\cite{Bi1}, was finally resolved by
the authors, the answer being positive if and only if
$n=1$~\cite{GoG}. 

In this paper, we study the braid groups of $\rp$, in particular the
splitting of the sequence~\reqref{split}, and the existence of a
section for the fibration $p$. These groups were first studied by Van
Buskirk~\cite{vB}, and more recently by Wang~\cite{Wa1}. Clearly
$P_1(\rp)=B_1(\rp)\cong \Z_2$. Van Buskirk showed that $P_2(\rp)$ is
isomorphic to the quaternion group $\mathcal{Q}_8$, $B_2(\rp)$ is a
generalised quaternion group of order~$16$, and for $n>2$, $P_n(\rp)$
and $B_n(\rp)$ are infinite. He also proved that these groups have
elements of finite order (including one of order~$2n$ in $B_n(\rp)$).
The torsion elements (although not their orders) of $B_n(\rp)$ were
characterised by Murasugi~\cite{M}. In~\cite{GG3}, we showed that for
$n\geq 2$, $B_n(\rp)$ has an element of order $\ell$ if and only if
$\ell$ divides $4n$ or $4(n-1)$, and that $P_n(\rp)$ has torsion
exactly $2$ and $4$. With respect to the splitting problem, Van
Buskirk showed that for all $n\geq 2$, neither the fibration
$\map{p}{F_n(\rp)}[F_1(\rp)]$ nor the homomorphism
$\map{p_{\ast}}{P_n(\rp)}[P_1(\rp)]$  admit a cross-section (for $p$,
this is a manifestation of the fixed point property of $\rp$), but
that the fibration $\map{p}{F_3(\rp)}[F_2(\rp)]$ admits a
cross-section, and hence so does the corresponding homomorphism. It
follows from~\reqref{split} that $P_3(\rp)$ is isomorphic to a
semi-direct product of $\pi_1(\rp \setminus \brak{x_1,x_2})$, which is
a free group $\F[2]$ of rank~$2$, by $P_2(\rp)$ which as we mentioned,
is isomorphic to $\mathcal{Q}_8$ (see~\cite{GG3} for an explicit
algebraic section). This fact will be used in the proof of
\repr{klein} (see \resec{presl}). Although there is no relation with
the braid groups of the sphere, it is a curious fact that the
commutator subgroup of $B_4(\St)$ is isomorphic to a semi-direct
product of $\mathcal{Q}_8$ by $\F[2]$~\cite{GG4}. In fact $B_n(\St)$
possesses subgroups isomorphic to $\mathcal{Q}_8$ if and only if
$n\geq 4$ is even~\cite{GG6}.

In~\cite{GG3}, we determined the homotopy type of the universal
covering space of $F_n(\rp)$. From this, we were able to deduce the
higher homotopy groups of $F_n(\rp)$. Using coincidence theory, we
then showed that for $n=2,3$ and $m\geq 4-n$, neither the fibration
nor the short exact sequence~\reqref{split} admit a section. More
precisely:
\begin{thm}[\cite{GG3}]\label{th:split23} Let $r\geq 4$ and $n=2,3$. Then: 
\begin{enumerate}[(a)]
\item the fibration $\map {p}{F_r(\rp)}[F_n(\rp)]$ does not admit a cross-section. 
\item the Fadell-Neuwirth pure braid group short exact sequence :
\begin{equation*} 
1 \to P_{r-n}(\rp\setminus\brak{x_1,\ldots,x_n})
\stackrel{i_{\ast}}{\to} P_r(\rp) \stackrel{p_{\ast}}{\to} P_n(\rp) \to 1
\end{equation*} 
does not split. 
\end{enumerate} 
\end{thm}

Apart from Van Buskirk's results for $F_n(\rp) \to F_1(\rp)$ and
$F_3(\rp) \to F_2(\rp)$ (published in~1966), no other results are
known concerning the splitting of~\reqref{split} for the pure braid
groups of $\rp$. The question is posed explicitly in the case 
$r=n+1$ on page~97 of \cite{vB}. In this paper, we give a complete answer. The main theorem is:
\begin{thm}\label{th:nosplit}
For all $n\geq 3$ and $m\geq 1$, the Fadell-Neuwirth pure braid group short exact sequence~\reqref{split}:
\begin{equation*}
1 \to P_m(\rp\setminus\brak{x_1, \ldots, x_n}) \to P_{n+m}(\rp) \stackrel{p_{\ast}}{\to} P_n(\rp) \to 1
\end{equation*}
does not split, and the fibration $\map{p}{F_{n+m}(\rp)}[F_n(\rp)]$ does not admit a section.
\end{thm}

Taking into account Van Buskirk's results and \reth{split23}, we deduce immediately the following corollary:
\begin{cor}
If $m,n\in \N$, the homomorphism $\map{p_{\ast}}{P_{n+m}(\rp)}[P_n(\rp)]$ and
the fibration $\map{p}{F_{n+m}(\rp)}[F_n(\rp)]$ admit a section if and only if
$n=2$ and $m=1$.\qed
\end{cor}
In other words, Van Buskirk's values ($n=2$ and $m=1$) are the only
ones for which a section exists (both on the geometric and the
algebraic level). The splitting problem for non-orientable surfaces
without boundary and of higher genus is the subject of work in
progress~\cite{GG5}. In the case of the Klein bottle, the existence of
a non-vanishing vector field implies that there always exists a
section, both geometric and algebraic (cf.~\cite{FaN}). 

This paper is organised as follows. In \resec{pres}, we start by
determining a presentation of $P_n(\rp)$ (\reth{basicpres}). To the
best of our knowledge, surprisingly this appears to be the first such
presentation in the literature (although Van Buskirk gave a
presentation of $B_n(\rp)$). 

In order to prove \reth{nosplit}, we argue by contradiction, and
suppose that there exists some $n\geq 3$ for which a section occurs.
As we indicate in \resec{proof}, it then suffices to study the case
$m=1$. The general strategy of the proof of \reth{nosplit} is based on
the following remark: if $H$ is any normal subgroup of $P_{n+1}(\rp)$
contained in $\ker{p_{\ast}}$, the quotiented short exact sequence \(
1\to \ker{p_{\ast}}/H \lhra P_{n+1}(\rp)/H \to P_n(\rp)\to 1 \) must also
split. In order to reach a contradiction, we seek such a subgroup
$H$ for which this short exact sequence does \emph{not} split. However
the choice of $H$ needed to achieve this is extremely delicate: if $H$
is too `small', the structure of the quotient $P_{n+1}(\rp)/H$ remains
complicated; on the other hand, if $H$ is too `big', we lose too much
information and cannot reach a conclusion. Taking a variety of
possible candidates for $H$, we observed in preliminary calculations
that the line between the two is somewhat fine. If $n$ is odd, we were able to show that the problem may be solved by taking the quotient $\ker{p_{\ast}}/H$ to be Abelianisation of $\ker{p_{\ast}}$ (which is a free Abelian group of rank $n$) modulo $2$, which is isomorphic to the direct sum of $n$ copies of $\Z_2$. However, this  insufficient for $n$ even.

With this in mind, in \resec{presl}, we study the quotient of
$P_{n+1}(\rp)$ by a certain normal subgroup $L$ which is contained in
$\ker{p_{\ast}}$ in the case $m=1$. A key step in the proof of
\reth{nosplit} is \repr{klein} where we show that $\ker{p_{\ast}}/L$
is isomorphic to $\Z^{n-1} \rtimes \Z$, the action being given by
multiplication by $-1$. This facilitates the calculations in
$P_{n+1}(\rp)/L$, whilst leaving just enough room for a contradiction.
This is accomplished in \resec{proof} where we show that the following
quotiented short exact sequence: 
\begin{equation*}
1 \to \ker{p_{\ast}}/L \to P_{n+1}(\rp)/L \to P_n(\rp) \to 1
\end{equation*}
does not split. 

\subsection*{Acknowledgements}

This work took place during the visit of the second author to the
Departmento de Matem\'atica do IME-Universidade de S\~ao Paulo during
the period 18\up{th}~June~--~18\up{th}~July 2006, and of the visit of
the first author to the Laboratoire de Math\'ematiques Emile Picard,
Université Paul Sabatier during the period
15\up{th}~November~--~16\up{th}~December 2006. This project was
supported by the international Cooperation USP/Cofecub project number
105/06.

\section{A presentation of $P_n(\rp)$}\label{sec:pres}

\begin{figure}[h]
\centering{\scalebox{0.6}{%
\input{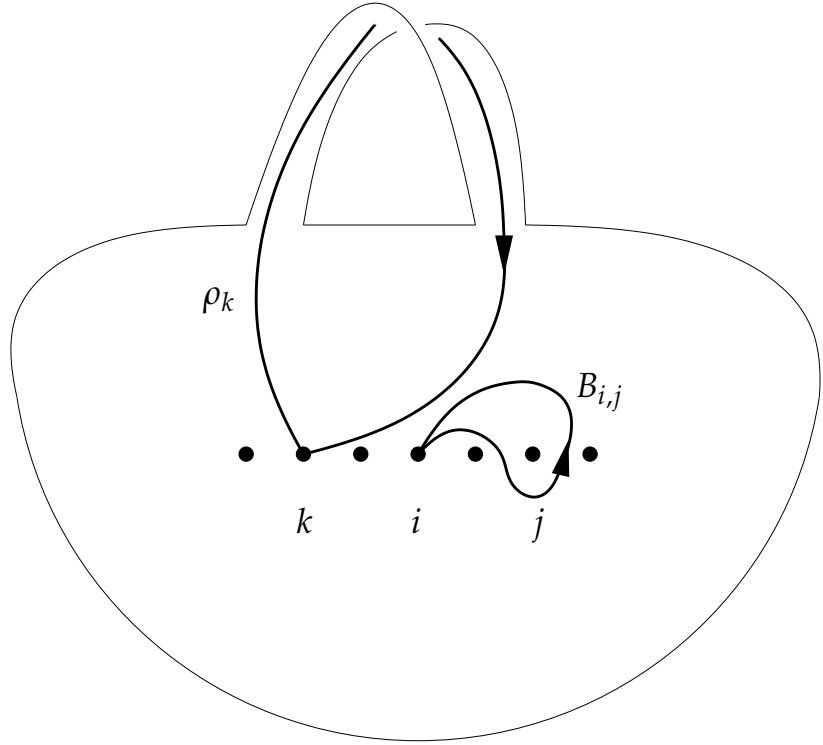}}} 
\caption{The generators $B_{i,j}$ and $\rho_k$ of $P_n(\rp)$.}
\label{fig:gens}
\end{figure}

If $n\in\N$ and $\dt\subseteq \rp$ is a topological disc, the
inclusion induces a (non-injective) homomorphism $\map
{\iota}{B_n(\dt)}[B_n(\rp)]$. If $\beta\in B_n(\dt)$ then we shall
denote its image $\iota(\beta)$ simply by $\beta$. For $1\leq i<j\leq
n$, we consider the following elements of $P_n(\rp)$:
\begin{equation*}\label{eq:genspn} 
B_{i,j}= \sigma_i^{-1}\cdots \sigma_{j-2}^{-1} \sigma_{j-1}^2 \sigma_{j-2} \cdots \sigma_i, 
\end{equation*} 
where $\sigma_1,\ldots, \sigma_{n-1}$ are the standard generators of
$B_n(\dt)$. The geometric braid corresponding to $B_{i,j}$ takes the
$i\up{th}$ string once around the $j\up{th}$ string in the positive
sense, with all other strings remaining vertical. For each $1\leq
k\leq n$, we define a generator $\rho_k$ which is represented
geometrically by a loop based at the $k\up{th}$ point and which goes
round the twisted handle. These elements are illustrated in
Figure~\ref{fig:gens} ($\rp$ minus a disc may be thought of as the
union of a disc and a twisted handle).

A presentation of $B_n(\rp)$ was first given by Van Buskirk in~\cite{vB}. Although presentations of braid groups of orientable and non-orientable surfaces have been the focus of several papers~\cite{Bi1,S,GM,Be}, we were not able to find an explicit presentation of $P_n(\rp)$ in the literature, so we derive one here. 

\begin{thm}\label{th:basicpres}
Let $n\in\N$. The following constitutes a presentation of pure braid group $P_n(\rp)$:
\begin{enumerate}
\item[\underline{\textbf{generators:}}] $B_{i,j}$, $1\leq i<j\leq n$, and $\rho_k$, $1\leq k\leq n$.
\item[\underline{\textbf{relations:}}]\mbox{}
\begin{enumerate}[(a)]
\item\label{it:rel1} the Artin relations between the $B_{i,j}$ emanating from those of $P_n(\dt)$:
\begin{equation*}
B_{r,s}B_{i,j}B_{r,s}^{-1}=
\begin{cases}
B_{i,j} & \text{if $i<r<s<j$ or $r<s<i<j$}\\
B_{i,j}^{-1} B_{r,j}^{-1}  B_{i,j} B_{r,j} B_{i,j} & \text{if $r<i=s<j$}\\
B_{s,j}^{-1} B_{i,j} B_{s,j} & \text{if $i=r<s<j$}\\
B_{s,j}^{-1}B_{r,j}^{-1} B_{s,j} B_{r,j} B_{i,j} B_{r,j}^{-1} B_{s,j}^{-1} B_{r,j} B_{s,j}   & \text{if $r<i<s<j$.}
\end{cases}
\end{equation*}
\item\label{it:rel2} for all $1\leq i<j\leq n$, $\rho_i\rho_j\rho_i^{-1} = \rho_j^{-1} B_{i,j}^{-1}  \rho_j^2$.
\item\label{it:rel3} for all $1\leq i\leq n$, the `surface relations' $\rho_i^2= B_{1,i}\cdots B_{i-1,i} B_{i,i+1} \cdots B_{i,n}$.
\item\label{it:rel4} for all $1\leq i<j\leq n$ and $1\leq k\leq n$, $k\neq j$,
\begin{equation*}\label{eq:relspn}
\rho_k B_{i,j}\rho_k^{-1}=
\begin{cases}
B_{i,j} & \text{if $j<k$ or $k<i$}\\
\rho_j^{-1} B_{i,j}^{-1} \rho_j & \text{if $k=i$}\\
\rho_j^{-1} B_{k,j}^{-1} \rho_j B_{k,j}^{-1} B_{i,j} B_{k,j} \rho_j^{-1} B_{k,j} \rho_j & \text{if $i<k<j$}.
\end{cases}
\end{equation*}
\end{enumerate}
\end{enumerate}
\end{thm}

\begin{proof}
We apply induction and standard results concerning the presentation of an extension (see Theorem~1, Chapter~13 of~\cite{J}).

First note that the given presentation is correct for $n=1$ ($P_1(\rp)= \pi_1(\rp)\cong \Z_2$), and $n=2$ ($P_2(\rp)\cong \mathcal{Q}_8$). So let $n\geq 2$, and suppose that $P_n(\rp)$ has the given presentation. Consider the corresponding Fadell-Neuwirth short exact sequence:
\begin{equation}\label{eq:fnm1}
1 \to \pi_1(\rp\setminus\brak{x_1, \ldots, x_n}) \to P_{n+1}(\rp) \stackrel{p_{\ast}}{\to} P_n(\rp) \to 1.
\end{equation}
In order to retain the symmetry of the presentation, we take the free group $\ker{p_{\ast}}$ to have the following one-relator presentation:
\begin{equation*}
\setangr{\rho_{n+1}, B_{1,n+1}, \ldots, B_{n,n+1}}{\rho_{n+1}^2= B_{1,n+1} \cdots B_{n,n+1}}.
\end{equation*}
Together with these generators of $\ker{p_{\ast}}$, the elements $B_{i,j}$, $1\leq i<j\leq n$, and $\rho_k$,  $1\leq k\leq n$, of $P_{n+1}(\rp)$ (which are coset representatives of the generators of $P_n(\rp)$)  form the required generating set of $P_{n+1}(\rp)$. 

There are three classes of relations of $P_{n+1}(\rp)$ which are
obtained as follows. The first consists of the single relation
$\rho_{n+1}^2= B_{1,n+1} \cdots B_{n,n+1}$ of $\ker{p_{\ast}}$. The
second class is obtained by rewriting the relators of the quotient in
terms of the coset representatives, and expressing the corresponding
element as a word in the generators of $\ker{p_{\ast}}$. In this way,
all of the relations of $P_n(\rp)$ lift directly to relations of
$P_{n+1}(\rp)$, with the exception of the surface relations which
become $\rho_i^2= B_{1,i}\cdots B_{i-1,i} B_{i,i+1} \cdots B_{i,n}
B_{i,n+1}$ for all $1\leq i\leq n$. Together with the relation of
$\ker{p_{\ast}}$, we obtain the complete set of surface relations
(relations~(\ref{it:rel3})) for $P_{n+1}(\rp)$.

The third class of relations is obtained by rewriting the conjugates
of the generators of $\ker{p_{\ast}}$ by the coset representatives in
terms of the generators of $\ker{p_{\ast}}$:
\begin{enumerate}[(i)]
\item\label{it:relex1} For all $1\leq i<j\leq n$ and $1\leq l\leq n$,
\begin{equation*}
B_{i,j} B_{l,n+1} B_{i,j}^{-1}= 
\begin{cases}
B_{l,n+1} & \text{if $l<i$ or $j<l$}\\
B_{l,n+1}^{-1} B_{i,n+1}^{-1}  B_{l,n+1} B_{i,n+1} B_{l,n+1} & \text{if $l=j$}\\
B_{j,n+1}^{-1} B_{l,n+1} B_{j,n+1} & \text{if $l=i$}\\
B_{j,n+1}^{-1}B_{i,n+1}^{-1} B_{j,n+1} B_{i,n+1} B_{l,n+1} B_{i,n+1}^{-1} B_{j,n+1}^{-1} B_{i,n+1} B_{j,n+1}   & \text{if $i<l<j$.}
\end{cases}
\end{equation*}
\item\label{it:relex2} $B_{i,j} \rho_{n+1} B_{i,j}^{-1}= \rho_{n+1}$ for all $1\leq i<j\leq n$.
\item\label{it:relex3} $\rho_k \rho_{n+1} \rho_k^{-1}= \rho_{n+1}^{-1} B_{k,n+1}^{-1}  \rho_{n+1}^2$ for all $1\leq k\leq n$.
\item\label{it:relex4} For all $1\leq k,l\leq n$,
\begin{equation*}
\rho_k B_{l,n+1}\rho_k^{-1}=
\begin{cases}
B_{l,n+1} & \text{if $k<l$}\\
\rho_{n+1}^{-1} B_{l,n+1}^{-1} \rho_{n+1} & \text{if $k=l$}\\
\rho_{n+1}^{-1} B_{k,n+1}^{-1} \rho_{n+1} B_{k,n+1}^{-1} B_{l,n+1} B_{k,n+1} \rho_{n+1}^{-1} B_{k,n+1} \rho_{n+1} & \text{if $l<k$.}
\end{cases}
\end{equation*}
\end{enumerate}
Then relations~(\ref{it:rel1}) for $P_{n+1}(\rp)$ are obtained from relations~(\ref{it:rel1}) for $P_{n}(\rp)$ and relations~(\ref{it:relex1}), relations~(\ref{it:rel2}) for $P_{n+1}(\rp)$ are obtained from relations~(\ref{it:rel2}) for $P_{n}(\rp)$ and relations~(\ref{it:relex3}), and relations~(\ref{it:rel4}) for $P_{n+1}(\rp)$ are obtained from relations~(\ref{it:rel4}) for $P_{n}(\rp)$, relations~(\ref{it:relex4}) and~(\ref{it:relex2}).
\end{proof}

For future use, it will be convenient at this point to record the following supplementary relations in $P_{n}(\rp)$ which are consequences of the presentation of \reth{basicpres}. Let $1\leq i<j \leq n$.
\begin{enumerate}[(I)]
\item\label{it:rel5} The action of the $\rho_i^{-1}$ on the $\rho_j$ may be deduced from that of $\rho_i$: $\rho_i^{-1}\rho_j\rho_i = B_{i,j}^{-1} \rho_j$.
\item By relations~(\ref{it:rel2}) and~(\ref{it:rel4}), we have:
\begin{align*}
\rho_i(B_{i,j}^{-1} \rho_j B_{i,j} \rho_j^{-1} B_{i,j}) \rho_i^{-1} &= \rho_j^{-1} B_{i,j} \rho_j \cdot \rho_j^{-1} B_{i,j}^{-1} \rho_j^2 \cdot \rho_j^{-1} B_{i,j}^{-1} \rho_j \cdot \rho_j^{-2}B_{i,j} \rho_j \cdot \rho_j^{-1} B_{i,j}^{-1} \rho_j
= B_{i,j}^{-1}.
\end{align*}
Hence $\rho_j B_{i,j} \rho_j^{-1}= B_{i,j} \rho_i^{-1} B_{i,j}^{-1} \rho_i B_{i,j}^{-1}$.
\item\label{it:rel6} From relations~(\ref{it:rel2}) and~(\ref{it:rel5}), we see that:
\begin{align*}
\rho_j \rho_i^{-1} \rho_j^{-1} &= \rho_i^{-1} \rho_j^{-1} B_{i,j}^{-1} \rho_j \rho_i \cdot \rho_i^{-1}= \rho_j^{-1} B_{i,j} \cdot \rho_i^{-1} B_{i,j}^{-1} \rho_i \cdot B_{i,j}^{-1} \rho_j \cdot \rho_i^{-1}= B_{i,j},
\end{align*}
so $\rho_j \rho_i \rho_j^{-1} = \rho_i B_{i,j}^{-1}$.
\item From relations~(\ref{it:rel5}) and~(\ref{it:rel4}), we obtain:
\begin{align*}
\rho_j^{-1} \rho_i \rho_j & = \rho_i \rho_j^{-1} B_{i,j} \rho_j  =\rho_i^2 B_{i,j}^{-1} \rho_i^{-1}.
\end{align*}
\end{enumerate}

\section{A presentation of the quotient $P_{n+1}(\rp)/L$}\label{sec:presl}

For $n\geq 2$, we have the Fadell-Neuwirth short exact sequence~(\ref{eq:fnm1}) whose kernel $K=\ker{p_{\ast}}$ is a free group of rank $n$ with basis $\rho_{n+1}, B_{1,n+1}, B_{2,n+1}, \ldots, B_{n-1,n+1}$. We first introduce a subgroup $L$ of $K$ which is normal in $P_{n+1}(\rp)$, from which we shall be able to prove Theorem~\ref{th:nosplit}. 

We define $L$ to be the normal closure in $P_{n+1}(\rp)$ of the following elements:
\begin{enumerate}[(i)]
\item $[B_{i,n+1}, B_{j,n+1}]$, where $1\leq i<j\leq n-1$, and
\item $[B_{i,n+1}, \rho_k]$, where $1\leq i\leq n-1$ and $1\leq k\leq n$.
\end{enumerate}
The elements $[B_{i,n+1}, B_{j,n+1}]$ clearly belong to $K$. The presentation of $P_n(\rp)$ given by \reth{basicpres} implies that:
\begin{equation*}
[B_{i,n+1}, \rho_k]=
\begin{cases}
1 & \text{if $k<i$}\\
B_{i,n+1}\rho_{n+1}^{-1} B_{i,n+1} \rho_{n+1} & \text{if $k=i$}\\
B_{i,n+1} \rho_{n+1}^{-1} B_{k,n+1}^{-1} \rho_{n+1} B_{k,n+1}^{-1} B_{i,n+1}^{-1} B_{k,n+1} \rho_{n+1}^{-1} B_{k,n+1} \rho_{n+1} & \text{if $i<k \leq n$.}
\end{cases} 
\end{equation*}
Thus $L$ is a (normal) subgroup of $K$.

Let $\map{g}{P_{n+1}(\rp)}[P_{n+1}(\rp)/L]$ denote the canonical
projection. For $i=1,\ldots, n-1$, let $A_i= g(B_{i,n+1})$. Apart from
these elements, if $x$ is a generator of $P_{n+1}(\rp)$, we shall not
distinguish notationally between $x$ and  $g(x)$. The quotient
$P_{n+1}(\rp)/L$ is generated by $\rho_1, \ldots, \rho_{n+1}$,
$B_{i,j}$, $1\leq i<j\leq n$, and $A_1,A_2, \ldots, A_{n-1}$ (we
delete $B_{n,n+1}$ from the list using the surface relation
$\rho_{n+1}^2= A_1 A_2 \cdots A_{n-1} B_{n,n+1}$, so
$B_{n,n+1}=A_{n-1}^{-1} \cdots A_2^{-1} A_1^{-1} \rho_{n+1}^2$). 

A presentation of $P_{n+1}(\rp)/L$ may obtained from that of $P_{n+1}(\rp)$ by adding the relations arising from the elements of $L$. We list those relations which are relevant for our description of  $P_{n+1}(\rp)/L$.
\begin{enumerate}[(a)]
\item The Artin relations between the $B_{i,j}$, $1\leq i<j\leq n$.
\item The relations of $P_n(\rp)$ between $\rho_i, \rho_j$, $1\leq i<j\leq n$.
\item The relations of $P_n(\rp)$ between $B_{i,j}$, $1\leq i<j\leq n$ and $\rho_k$, $1\leq k\leq n$.

The following two sets of relations arise from the definition of $L$:

\item\label{it:comma} $A_i\comm A_j$, $1\leq i <j \leq n-1$ (the symbol $\comm$ is used to mean that the given elements commute).
\item $A_i\comm \rho_j$, $i=1,\ldots, n-1$ and $j=1,\ldots ,n$.

\item The surface relations:
\begin{align*}
\rho_i^2 & = B_{1,i}\cdots B_{i-1,i} B_{i,i+1} \cdots B_{i,n} A_i \;\text{for $i=1,\ldots, n-1$}\\
\rho_n^2 &= B_{1,n}B_{2,n} \cdots B_{n-1,n}\cdot A_{n-1}^{-1}\cdots A_2^{-1}A_1^{-1}\rho_{n+1}^2.
\end{align*}

\item\label{it:klein} For $i=1,\ldots, n-1$, $\rho_{n+1} A_i \rho_{n+1}^{-1}= A_i^{-1}$ (since $\rho_{n+1}^{-1} A_i^{-1} \rho_{n+1}= \rho_i A_i\rho_i^{-1}=A_i$).

\item[] The following relations are implied by the above relations:
\begin{enumerate}[--]
\item for $1\leq i<j\leq n-1$, $\rho_j A_i \rho_j^{-1} = \rho_{n+1}^{-1} A_j^{-1} \rho_{n+1} A_j^{-1} A_i A_j \rho_{n+1}^{-1} A_j \rho_{n+1}$ (both are equal to $A_i$).
\item for $1\leq i\leq n-1$, $\rho_{n+1} A_i\rho_{n+1}^{-1}= A_i \rho_i^{-1} A_i^{-1} \rho_i A_i^{-1}$ (both are equal to $A_i^{-1}$)
\end{enumerate}

\item For $1\leq j\leq n-1$,
\begin{equation*}
\rho_i \rho_{n+1} \rho_i^{-1}= \rho_{n+1}^{-1}
A_i^{-1}\rho_{n+1}^2=A_i \rho_{n+1}=\rho_{n+1} A_i^{-1}.
\end{equation*}
From these relations, it follows that $\rho_i\comm \rho_{n+1}^2$ for $i=1,\ldots,n-1$.

\item $\rho_n \rho_{n+1} \rho_n^{-1}  = \rho_{n+1}^{-1} B_{n,n+1}^{-1} \rho_{n+1}^2 =\rho_{n+1}^{-1} \rho_{n+1}^{-2} A_1 \cdots A_{n-1} \rho_{n+1}^2=
A_1^{-1} \cdots A_{n-1}^{-1} \rho_{n+1}^{-1}$.
From this relation, it follows that $\rho_n\rho_{n+1}^2 \rho_n^{-1}= \rho_{n+1}^{-2}$.

\item[] For $i=1\ldots, n-1$, the following relations are implied by the above relations:
\begin{equation*}
\text{$\rho_n A_i \rho_n^{-1} = \rho_{n+1}^{-1} B_{n,n+1}^{-1} \rho_{n+1} B_{n,n+1}^{-1} A_i B_{n,n+1} \rho_{n+1}^{-1} B_{n,n+1} \rho_{n+1}$ (both are equal to $A_i$).}
\end{equation*}
\end{enumerate}

\begin{prop}\label{prop:klein}
The quotient group $K/L$ has a presentation of the form:
\begin{enumerate}
\item[\textbf{generators:}] $A_1,\ldots, A_{n-1},\rho_{n+1}$.
\item[\textbf{relations:}] $A_i \comm A_j$ for $1\leq i<j \leq n-1$, and $\rho_{n+1} A_i \rho_{n+1}^{-1}= A_i^{-1}$ for $1\leq i \leq n-1$.
\end{enumerate}
In particular, $K/L$ is isomorphic to $\Z^{n-1} \rtimes \Z$, the action being given by multiplication by $-1$.
\end{prop}

Hence the other relations in $P_{n+1}(\rp)/L$ (which involve only elements from $P_n(\rp)$) do not add any further relations to the quotient $K/L$.

\begin{proof}[Proof of Proposition~\ref{prop:klein}]
Clearly $A_1,\ldots, A_{n-1},\rho_{n+1}$ generate $K/L$, and from  relations~(\ref{it:comma}) and~(\ref{it:klein}) of $P_{n+1}(\rp)/L$, they are subject to the given relations. Consider the following commutative diagram of short exact sequences:
\begin{equation}\label{eq:seskl}
\begin{xy}*!C\xybox{\xymatrix{%
1\ar[r] & K \txt{~~}\ar@{^{(}->}[r] \ar^{g\left\lvert_K \right.}[d] & P_{n+1}(\rp) \ar^{p_{\ast}}[r] \ar^{g}[d] & P_n(\rp) \ar[r] \ar@{=}[d] & 1\\
1\ar[r] & K/L \txt{~~} \ar^(0.35){{\iota}}@{^{(}->}[r]  & P_{n+1}(\rp)/L \ar^(0.55){\overline{p}_{\ast}}[r] & P_n(\rp) \ar[r]  & 1,}}
\end{xy}
\end{equation}
where $\iota$ is the inclusion of $K/L$ in $P_{n+1}(\rp)/L$, and $\overline{p}_{\ast}$ is the homomorphism induced by $p_{\ast}$. Let $\Gamma$ be the group with presentation:
\begin{equation*}
\Gamma=\setangr{\alpha_1, \ldots, \alpha_{n-1}, \rho}{\text{$\alpha_i \comm \alpha_j$ for $1\leq i<j \leq n-1$, and $\rho \alpha_i \rho^{-1}= \alpha_i^{-1}$}}.
\end{equation*}
So $\Gamma$ is isomorphic to $\Z^{n-1} \rtimes \Z$, where the action is given by multiplication by $-1$. The map $\map{f}{\Gamma}[K/L]$ defined on the generators of $\Gamma$ by $f(\alpha_i)=A_i$ for $i=1, \ldots ,n-1$, and $f(\rho)=\rho_{n+1}$, extends to a surjective homomorphism. We claim that $f$ is an isomorphism, which will prove the proposition. To prove the claim, it suffices to show that $\iota \circ f$ is injective. Let $w\in \ker{\iota \circ f}$. Then we may write $w$ uniquely in the form $w= \rho^{m_0} \alpha_1^{m_1} \cdots \alpha_{n-1}^{m_{n-1}}$, where $m_0,m_1, \ldots, m_{n-1}\in \Z$, and so $\iota \circ f(w)=\rho_{n+1}^{m_0} A_1^{m_1} \cdots A_{n-1}^{m_{n-1}}=1$ in $P_{n+1}(\rp)/L$. 

Let $z=\rho_{n+1}^{m_0} B_{1,{n+1}}^{m_1} \cdots B_{n-1,n+1}^{m_{n-1}} \in P_{n+1}(\rp)$. Since $g(z)=\iota \circ f(w)=1$, we must have $z\in L$. Now $L$ is the normal closure in $P_{n+1}(\rp)$ of the following elements:
\begin{enumerate}[--]
\item $c_{j,k}=[B_{j,n+1}, B_{k,n+1}]$, where $1\leq j<k\leq n-1$, 
\item $d_j= [B_{j,n+1}, \rho_j]= B_{j,n+1} \rho_{n+1}^{-1} B_{j,n+1} \rho_{n+1}$, where $1\leq j\leq n-1$,
\item $e_{j,k}= [B_{j,n+1}, \rho_k]= B_{j,n+1} \rho_{n+1}^{-1} B_{k,n+1}^{-1} \rho_{n+1} B_{k,n+1}^{-1} B_{j,n+1}^{-1} B_{k,n+1} \rho_{n+1}^{-1} B_{k,n+1} \rho_{n+1}$, where $1\leq j<k\leq n$. 
\end{enumerate}
Hence $z$ may be written as a product of conjugates of the elements $c_{j,k}$, $d_j$, $e_{j,k}$, and their inverses.

For $i=1,\ldots, n-1$, let $\map{\pi_i}{P_{n+1}(\rp)}[P_3(\rp,
(p_i,p_n,p_{n+1}))]$ be the projection obtained geometrically by
forgetting all of the strings, with the exception of the $i\up{th}$,
$n\up{th}$ and $(n+1)\up{st}$ strings (here $P_3(\rp,
(p_i,p_n,p_{n+1}))$ denotes the fundamental group of $F_3(\rp)$ taking
the basepoint to be $(p_i,p_n,p_{n+1})$). We interpret $P_3(\rp,
(p_i,p_n,p_{n+1}))$ as the semi-direct product $\F[2](B_{i,n+1},
\rho_{n+1}) \rtimes P_2(\rp, (p_i,p_n))$~\cite{vB}. Under $\pi_i$, the
elements $c_{j,k},d_j, e_{j,k}$ (for the allowed values of $j$ and
$k$) are all sent to the trivial element, with the exception of the
two elements $d_i$ and $e_{i,n}$. Set $h_i=\pi_i(d_i)= B_{i,n+1}
\rho_{n+1}^{-1} B_{i,n+1} \rho_{n+1}\in \F[2](B_{i,n+1}, \rho_{n+1})$.
Since $B_{i,n+1} B_{n,n+1}= \rho_{n+1}^2$ in $P_3(\rp,
(p_i,p_n,p_{n+1}))$, we have $B_{n,n+1}= B_{i,n+1}^{-1} \rho_{n+1}^2$.
Hence:
\begin{align*}
\pi_i(e_{i,n}) =& B_{i,n+1} \rho_{n+1}^{-1} B_{n,n+1}^{-1} \rho_{n+1} B_{n,n+1}^{-1} B_{i,n+1}^{-1} B_{n,n+1} \rho_{n+1}^{-1} B_{n,n+1} \rho_{n+1}\\
=& B_{i,n+1} \rho_{n+1}^{-1} \rho_{n+1}^{-2} B_{i,n+1} \rho_{n+1} \rho_{n+1}^{-2} B_{i,n+1} B_{i,n+1}^{-1} B_{i,n+1}^{-1} \rho_{n+1}^2 \rho_{n+1}^{-1} B_{i,n+1}^{-1} \rho_{n+1}^2 \rho_{n+1}\\
=& B_{i,n+1} \rho_{n+1}^{-1} B_{i,n+1} \rho_{n+1} \cdot \rho_{n+1}^{-1} B_{i,n+1}^{-1} \rho_{n+1}^{-1} B_{i,n+1}^{-1} \rho_{n+1} B_{i,n+1}^{-1} B_{i,n+1} \rho_{n+1} \cdot\\
& \rho_{n+1}^{-2} B_{i,n+1} \rho_{n+1}^{-1} B_{i,n+1} \rho_{n+1} \rho_{n+1}^2 \cdot \rho_{n+1}^{-3} \rho_{n+1}^{-1} B_{i,n+1}^{-1} \rho_{n+1} B_{i,n+1}^{-1} \rho_{n+1}^3\\
=& h_i\cdot \rho_{n+1}^{-1} B_{i,n+1}^{-1} h_i^{-1} B_{i,n+1} \rho_{n+1} \cdot \rho_{n+1}^{-2} h_i \rho_{n+1}^2 \cdot \rho_{n+1}^{-3} h_i^{-1} \rho_{n+1}^3.
\end{align*}
Thus $\pi_i(z)$ may be written as a product of conjugates in $P_3(\rp, (p_i,p_n,p_{n+1}))$ of $h_i^{\pm 1}$:
\begin{equation}\label{eq:pi1z}
\pi_i(z)= \rho_{n+1}^{m_0} B_{i,{n+1}}^{m_i}= \prod_{j=1}^l \; w_j h_i^{\mu(j)} w_j^{-1},
\end{equation}
where $l\in\N$, $w_j\in P_3(\rp, (p_i,p_n,p_{n+1}))$, and $\mu(j)\in \brak{1,-1}$. We claim that each $w_j h_i^{\mu(j)} w_j^{-1}$ is in fact a conjugate in $\F[2](B_{i,n+1}, \rho_{n+1})$ of $h_i^{\pm 1}$. This follows by studying the action of the generators $\rho_i$ and $\rho_n$ of $P_2(\rp, (p_i,p_n))$ on the basis of $\F[2](B_{i,n+1}, \rho_{n+1})$:
\begin{align*}
\rho_i h_i \rho_i^{-1}=& \rho_i B_{i,n+1} \rho_{n+1}^{-1} B_{i,n+1} \rho_{n+1} \rho_i^{-1}\\
=& \rho_{n+1}^{-1} B_{i,n+1}^{-1} \rho_{n+1} \cdot \rho_{n+1}^{-2} B_{i,n+1} \rho_{n+1} \cdot \rho_{n+1}^{-1} B_{i,n+1}^{-1} \rho_{n+1} \cdot \rho_{n+1}^{-1} B_{i,n+1}^{-1}  \rho_{n+1}^2\\
=& \rho_{n+1}^{-1} B_{i,n+1}^{-1} \rho_{n+1}^{-1}  B_{i,n+1}^{-1} \rho_{n+1}^2\\
=& \rho_{n+1}^{-1} B_{i,n+1}^{-1} \cdot \rho_{n+1}^{-1}  B_{i,n+1}^{-1} \rho_{n+1} B_{i,n+1}^{-1} \cdot B_{i,n+1} \rho_{n+1} = \rho_{n+1}^{-1} B_{i,n+1}^{-1} h_i^{-1} B_{i,n+1} \rho_{n+1},
\end{align*}
and
\begin{align*}
\rho_n h_i \rho_n^{-1}=& \rho_n B_{i,n+1} \rho_{n+1}^{-1} B_{i,n+1} \rho_{n+1} \rho_n^{-1}\\
=& \rho_{n+1}^{-1} B_{n,n+1}^{-1} \rho_{n+1} B_{n,n+1}^{-1} B_{i,n+1} B_{n,n+1} \rho_{n+1}^{-1} B_{n,n+1} \rho_{n+1} \cdot \rho_{n+1}^{-2} B_{n,n+1} \rho_{n+1} \cdot \\
& \rho_{n+1}^{-1} B_{n,n+1}^{-1} \rho_{n+1} B_{n,n+1}^{-1} B_{i,n+1} B_{n,n+1} \rho_{n+1}^{-1} B_{n,n+1} \rho_{n+1} \cdot \rho_{n+1}^{-1} B_{n,n+1}^{-1}  \rho_{n+1}^2\\
=& \rho_{n+1}^{-1} B_{n,n+1}^{-1} \rho_{n+1} B_{n,n+1}^{-1} B_{i,n+1} B_{n,n+1} \rho_{n+1}^{-1} B_{i,n+1} B_{n,n+1} \rho_{n+1}\\
=& \rho_{n+1}^{-1} \rho_{n+1}^{-2} B_{i,n+1} \rho_{n+1} \rho_{n+1}^{-2} B_{i,n+1} B_{i,n+1} B_{i,n+1}^{-1} \rho_{n+1}^2 \rho_{n+1}^{-1} B_{i,n+1} B_{i,n+1}^{-1} \rho_{n+1}^2 \rho_{n+1}\\
=& \rho_{n+1}^{-3} B_{i,n+1} \rho_{n+1}^{-1} B_{i,n+1} \rho_{n+1} \rho_{n+1}^3 = \rho_{n+1}^{-3} h_i \rho_{n+1}^3,
\end{align*} 
again using the fact that $B_{i,n+1}B_{n,n+1}=\rho_{n+1}^2$ in $P_3(\rp, (p_i,p_n,p_{n+1}))$. Thus the $w_j$ of equation~(\ref{eq:pi1z}) may be taken as belonging to $\F[2](B_{i,n+1}, \rho_{n+1})$. We now project $\F[2](B_{i,n+1}, \rho_{n+1})$ onto the Klein bottle group $\setangr{B_{i,n+1}, \rho_{n+1}}{\rho_{n+1}^{-1} B_{i,n+1} \rho_{n+1}=B_{i,n+1}^{-1}}$ in the obvious manner. Since $h_i$ belongs to the kernel of this projection, the right hand-side of equation~(\ref{eq:pi1z}) is sent to the trivial element, while the left hand-side is sent to $\rho_{n+1}^{m_0} B_{i,n+1}^{m_i}$. It follows that $m_0=m_i=0$ for all $i=1,\ldots, n-1$. This proves the injectivity of $\iota \circ f$, and so completes the proof of the proposition.
\end{proof}

\section{Proof of Theorem~\ref{th:nosplit}}\label{sec:proof}

We are now ready to give the proof of the main theorem of the paper.

\begin{proof}[Proof of Theorem~\ref{th:nosplit}]
Let $n\geq 3$. For $m\geq 1$, let $\map{p_{\ast}^{(m)}}{P_{n+m}(\rp)}[P_n(\rp)]$ denote the usual projection. Suppose first that $m\geq 2$, and consider the following commutative diagram of short exact sequences:
\begin{equation*}
\xymatrix{%
1\ar[r] & P_m(\rp \setminus \brak{x_1,\ldots, x_n}) \txt{~~}\ar@{^{(}->}[r] \ar^{\psi\left\lvert_{P_m\left(\rp \setminus \brak{x_1,\ldots, x_n}\right)}\right.}[d] & P_{n+m}(\rp) \ar^(0.55){p_{\ast}^{(m)}}[r] \ar^{\psi}[d] & P_n(\rp) \ar[r] \ar@{=}[d] & 1\\
1\ar[r] & P_1(\rp \setminus \brak{x_1,\ldots, x_n}) \txt{~~} \ar@{^{(}->}[r]  & P_{n+1}(\rp) \ar^(0.55){p_{\ast}^{(1)}}[r] & P_n(\rp) \ar[r]  & 1,}
\end{equation*}
where $\psi$ is the homomorphism which forgets the last $m-1$ strings. If $p_{\ast}^{(m)}$ admits a section $s_{\ast}^{(m)}$ then $\psi \circ s_{\ast}^{(m)}$ is a section for $p_{\ast}^{(1)}$. In other words, if the upper short exact sequence splits then so does the lower one.

Since we shall be arguing for a contradiction, we are reduced to considering the case $m=1$. Set $p_{\ast}= p_{\ast}^{(1)}$, and suppose that $p_{\ast}$ admits a section which we shall denote by $s_{\ast}$. Consider the short exact sequence~(\ref{eq:seskl}). Since $p_{\ast}$ admits a section then so does $\overline{p}_{\ast}$; we denote its section by $\overline{s}_{\ast}$. So $\overline{p}_{\ast}(\rho_i)= \rho_i$ for $i=1,\ldots, n$, and $\overline{p}_{\ast}(B_{i,j})= B_{i,j}$ for $1\leq i<j\leq n$ (recall that we do not distinguish notationally between the generators of $P_{n+1}(\rp)/L$ and the corresponding generators of $P_n(\rp)$). Thus we obtain:
\begin{equation}\label{eq:defcoef}
\left.
\begin{aligned}
\overline{s}_{\ast}(\rho_i) &= \text{$\rho_{n+1}^{\indalp{i}{0}} A_1^{\indalp{i}{1}} \cdots A_{n-1}^{\indalp{i}{n-1}} \cdot \rho_i$ for $i=1,\ldots, n$}\\
\overline{s}_{\ast}(B_{i,j}) &= \text{$\rho_{n+1}^{\indbet{i}{j}{0}} A_1^{\indbet{i}{j}{1}} \cdots A_{n-1}^{\indbet{i}{j}{n-1}} \cdot B_{i,j}$ for $1\leq i <j \leq n$,}
\end{aligned}
\right\}
\end{equation}
where $\indalp{i}{k}, \indbet{i}{j}{k}\in \Z$. For $x\in \Z$, set
\begin{align*}
& \epsilon(x)=
\begin{cases}
1 & \text{if $x$ is even}\\
-1 & \text{if $x$ is odd,}
\end{cases} && \text{and} && \delta(x)=
\begin{cases}
0 & \text{if $x$ is even}\\
-1 & \text{if $x$ is odd.}
\end{cases}
\end{align*}
Then $\epsilon(x)=2\delta(x)+1$, $\epsilon(x)\delta(x)= -\delta(x)$, $\delta(x)= \delta(-x)$, $\epsilon(x)=\epsilon(-x)$ and for $i=1,\ldots, n-1$ and $k\in \Z$, we have: 
\begin{align*}
\rho_{n+1}^k A_i \rho_{n+1}^{-k} &= A_i^{\epsilon(k)}\\
\rho_i \rho_{n+1}^k \rho_i^{-1} &= \rho_{n+1}^k A_i^{\delta(k)}\\
\rho_i^{-1} \rho_{n+1}^k \rho_i &= \rho_{n+1}^k A_i^{-\delta(k)}\\
\rho_n \rho_{n+1}^k \rho_n^{-1} &= \rho_{n+1}^{-k} A_1^{-\delta(k)} \cdots A_{n-1}^{-\delta(k)}\\
\rho_n^{-1} \rho_{n+1}^k \rho_n &= \rho_{n+1}^{-k} A_1^{\delta(k)} \cdots A_{n-1}^{\delta(k)},
\end{align*}
using the relations of $P_{n+1}(\rp)/L$ given in \resec{presl}.

We now calculate the images in $P_{n+1}(\rp)/L$ by $\overline{s}_{\ast}$ of the following relations of $P_n(\rp)$. This will allow us to obtain information about the coefficients defined in \req{defcoef}.
\begin{enumerate}[(a)]
\item We start with the relation \underline{$\rho_j\rho_i\rho_j^{-1}=\rho_i B_{i,j}^{-1}$} in $P_n(\rp)$, where $1\leq i<j\leq n-1$.
\begin{align*}
\overline{s}_{\ast}(\rho_i B_{i,j}^{-1}) =& \rho_{n+1}^{\indalp{i}{0}} A_1^{\indalp{i}{1}} \cdots A_{n-1}^{\indalp{i}{n-1}} \rho_i \cdot
B_{i,j}^{-1} A_{n-1}^{-\indbet{i}{j}{n-1}} \cdots A_1^{-\indbet{i}{j}{1}}  \rho_{n+1}^{-\indbet{i}{j}{0}}\\
=& \rho_{n+1}^{\indalp{i}{0}} A_1^{\indalp{i}{1}} \cdots A_{n-1}^{\indalp{i}{n-1}} \rho_i \rho_{n+1}^{-\indbet{i}{j}{0}}
B_{i,j}^{-1} A_{n-1}^{-\epsilon(\indbet{i}{j}{0}) \indbet{i}{j}{n-1}} \cdots A_1^{-\epsilon(\indbet{i}{j}{0}) \indbet{i}{j}{1}}\\
=& \rho_{n+1}^{\indalp{i}{0}} A_1^{\indalp{i}{1}} \cdots A_{n-1}^{\indalp{i}{n-1}}  \rho_{n+1}^{-\indbet{i}{j}{0}} A_i^{\delta(\indbet{i}{j}{0})}
\rho_i  A_{n-1}^{-\epsilon(\indbet{i}{j}{0}) \indbet{i}{j}{n-1}} \cdots A_1^{-\epsilon(\indbet{i}{j}{0}) \indbet{i}{j}{1}} B_{i,j}^{-1}\\
=& \rho_{n+1}^{\indalp{i}{0}-\indbet{i}{j}{0}} A_1^{\epsilon(\indbet{i}{j}{0}) \indalp{i}{1}} \cdots A_{n-1}^{\epsilon(\indbet{i}{j}{0}) \indalp{i}{n-1}}  A_i^{\delta(\indbet{i}{j}{0})}
\cdot\\
&  A_{n-1}^{-\epsilon(\indbet{i}{j}{0}) \indbet{i}{j}{n-1}} \cdots A_1^{-\epsilon(\indbet{i}{j}{0}) \indbet{i}{j}{1}} \rho_i B_{i,j}^{-1}\\
=& \rho_{n+1}^{\indalp{i}{0}-\indbet{i}{j}{0}} 
A_1^{\epsilon(\indbet{i}{j}{0}) (\indalp{i}{1} - \indbet{i}{j}{1})} \cdots
A_i^{\epsilon(\indbet{i}{j}{0}) (\indalp{i}{i} - \indbet{i}{j}{i}) +\delta(\indbet{i}{j}{0})}\cdot\\
& \cdots A_{n-1}^{\epsilon(\indbet{i}{j}{0}) (\indalp{i}{n-1}-\indbet{i}{j}{n-1})} \rho_i B_{i,j}^{-1}.
\end{align*}

\begin{align*}
\overline{s}_{\ast}(\rho_j\rho_i\rho_j^{-1}) =& 
\rho_{n+1}^{\indalp{j}{0}} A_1^{\indalp{j}{1}} \cdots A_{n-1}^{\indalp{j}{n-1}} \rho_j \cdot
\rho_{n+1}^{\indalp{i}{0}} A_1^{\indalp{i}{1}} \cdots A_{n-1}^{\indalp{i}{n-1}} \cdot \rho_i \cdot
\rho_j^{-1} A_{n-1}^{-\indalp{j}{n-1}} \cdots A_1^{-\indalp{j}{1}} \rho_{n+1}^{-\indalp{j}{0}}\\
=&
\rho_{n+1}^{\indalp{j}{0}+ \indalp{i}{0}} A_1^{\epsilon(\indalp{i}{0}) \indalp{j}{1}} \cdots A_{n-1}^{\epsilon(\indalp{i}{0}) \indalp{j}{n-1}}
A_j^{\delta(\indalp{i}{0})} \rho_j A_1^{\indalp{i}{1}} \cdots A_{n-1}^{\indalp{i}{n-1}} \rho_{n+1}^{-\indalp{j}{0}} A_i^{\delta(\indalp{j}{0})} \rho_i \cdot\\
& A_j^{-\delta(\indalp{j}{0})} \rho_j^{-1} A_{n-1}^{-\epsilon(\indalp{j}{0}) \indalp{j}{n-1}} \cdots A_1^{-\epsilon(\indalp{j}{0}) \indalp{j}{1}}\\
=&
\rho_{n+1}^{\indalp{i}{0}} A_1^{\epsilon(\indalp{j}{0}) \epsilon(\indalp{i}{0}) \indalp{j}{1}} \cdots A_{n-1}^{\epsilon(\indalp{j}{0}) \epsilon(\indalp{i}{0}) \indalp{j}{n-1}}
A_j^{\epsilon(\indalp{j}{0}) \delta(\indalp{i}{0})} 
A_j^{\delta(\indalp{j}{0})}\cdot\\
& A_1^{\epsilon(\indalp{j}{0}) \indalp{i}{1}} \cdots A_{n-1}^{\epsilon(\indalp{j}{0}) \indalp{i}{n-1}}  A_i^{\delta(\indalp{j}{0})}  
A_j^{-\delta(\indalp{j}{0})}  A_{n-1}^{-\epsilon(\indalp{j}{0}) \indalp{j}{n-1}} \cdots A_1^{-\epsilon(\indalp{j}{0}) \indalp{j}{1}} \rho_j \rho_i \rho_j^{-1}\\
=&
\rho_{n+1}^{\indalp{i}{0}} 
A_1^{\epsilon(\indalp{j}{0}) \left(\indalp{j}{1}(\epsilon(\indalp{i}{0})-1) + \indalp{i}{1}\right)} \cdots
A_i^{\epsilon(\indalp{j}{0}) \left(\indalp{j}{i}(\epsilon(\indalp{i}{0})-1) + \indalp{i}{i}\right)+ \delta(\indalp{j}{0})}\cdot\\
& \cdots 
A_j^{\epsilon(\indalp{j}{0}) \left(\indalp{j}{j}(\epsilon(\indalp{i}{0})-1) + \indalp{i}{j}+\delta(\indalp{i}{0})\right)} \cdots 
A_{n-1}^{\epsilon(\indalp{j}{0}) \left( \indalp{j}{n-1}(\epsilon(\indalp{i}{0})-1) + \indalp{i}{n-1}\right)}
\rho_j \rho_i \rho_j^{-1}.
\end{align*}

Comparing coefficients in $K/L$, we obtain:
\begin{align}
& \text{$\indbet{i}{j}{0}=0$, so $\epsilon(\indbet{i}{j}{0})=1$ and $\delta(\indbet{i}{j}{0})=0$ for all $1\leq i<j\leq n-1$}\label{eq:rho2rho1}\\
& \text{$\epsilon(\indalp{j}{0}) \indalp{j}{k}
(\epsilon(\indalp{i}{0})-1) + \indalp{i}{k}(\epsilon(\indalp{j}{0})-1)= -\indbet{i}{j}{k}$ for all $k=1,\ldots,n-1$, $k\neq i,j$}\notag\\
& \epsilon(\indalp{j}{0}) \indalp{j}{i}(\epsilon(\indalp{i}{0})-1) + \indalp{i}{i}(\epsilon(\indalp{j}{0})-1)+ \delta(\indalp{j}{0})= -\indbet{i}{j}{i}\notag\\
& \epsilon(\indalp{j}{0}) \indalp{j}{j}(\epsilon(\indalp{i}{0})-1) + \indalp{i}{j}(\epsilon(\indalp{j}{0})-1)+\epsilon(\indalp{j}{0}) \delta(\indalp{i}{0})= -\indbet{i}{j}{j}.\notag
\end{align}
In particular, the coefficient $\indbet{i}{j}{0}$ of $\rho_{n+1}$ in $\overline{s}_{\ast}(B_{i,j})$ is zero. Also, since $\epsilon(x)-1$ is even for all $x\in\Z$, $\indbet{i}{j}{k}\equiv 0 \pmod 2$ for all $k\neq i,j$, and for all $1\leq i<j\leq n-1$,
\begin{align}
\indbet{i}{j}{i} &\equiv \delta(\indalp{j}{0})  \pmod 2\label{eq:iji}\\
\indbet{i}{j}{j} &\equiv \delta(\indalp{i}{0}) \pmod 2.\label{eq:ijj}
\end{align}

\item Now consider the relation \underline{$\rho_n\rho_i\rho_n^{-1}=\rho_i B_{i,n}^{-1}$} in $P_n(\rp)$, where $1\leq i\leq n-1$.
\begin{align*}
\overline{s}_{\ast}(\rho_i B_{i,n}^{-1}) =& \rho_{n+1}^{\indalp{i}{0}} A_1^{\indalp{i}{1}} \cdots A_{n-1}^{\indalp{i}{n-1}} \rho_i \cdot
B_{i,n}^{-1} A_{n-1}^{-\indbet{i}{n}{n-1}} \cdots A_1^{-\indbet{i}{n}{1}}  \rho_{n+1}^{-\indbet{i}{n}{0}}\\
=& \rho_{n+1}^{\indalp{i}{0}} A_1^{\indalp{i}{1}} \cdots A_{n-1}^{\indalp{i}{n-1}} \rho_i \rho_{n+1}^{-\indbet{i}{n}{0}}
B_{i,n}^{-1} A_{n-1}^{-\epsilon(\indbet{i}{n}{0}) \indbet{i}{n}{n-1}} \cdots A_1^{-\epsilon(\indbet{i}{n}{0}) \indbet{i}{n}{1}}\\
=& \rho_{n+1}^{\indalp{i}{0}} A_1^{\indalp{i}{1}} \cdots A_{n-1}^{\indalp{i}{n-1}}  \rho_{n+1}^{-\indbet{i}{n}{0}} A_i^{\delta(\indbet{i}{n}{0})}
\rho_i  A_{n-1}^{-\epsilon(\indbet{i}{n}{0}) \indbet{i}{n}{n-1}} \cdots A_1^{-\epsilon(\indbet{i}{n}{0}) \indbet{i}{n}{1}} B_{i,n}^{-1}\\
=& \rho_{n+1}^{\indalp{i}{0}-\indbet{i}{n}{0}} A_1^{\epsilon(\indbet{i}{n}{0}) \indalp{i}{1}} \cdots A_{n-1}^{\epsilon(\indbet{i}{n}{0}) \indalp{i}{n-1}}  A_i^{\delta(\indbet{i}{n}{0})}\cdot\\
&  A_{n-1}^{-\epsilon(\indbet{i}{n}{0}) \indbet{i}{n}{n-1}} \cdots A_1^{-\epsilon(\indbet{i}{n}{0}) \indbet{i}{n}{1}} \rho_i B_{i,n}^{-1}\\
=& \rho_{n+1}^{\indalp{i}{0}-\indbet{i}{n}{0}} 
A_1^{\epsilon(\indbet{i}{n}{0}) (\indalp{i}{1} - \indbet{i}{n}{1})} \cdots
A_i^{\epsilon(\indbet{i}{n}{0}) (\indalp{i}{i} - \indbet{i}{n}{i}) +\delta(\indbet{i}{n}{0})} \cdot\\
& \cdots A_{n-1}^{\epsilon(\indbet{i}{n}{0}) (\indalp{i}{n-1}-\indbet{i}{n}{n-1})} \rho_i B_{i,n}^{-1}.
\end{align*}

\begin{align*}
\overline{s}_{\ast}(\rho_n\rho_i\rho_n^{-1}) =& 
\rho_{n+1}^{\indalp{n}{0}} A_1^{\indalp{n}{1}} \cdots A_{n-1}^{\indalp{n}{n-1}} \rho_n \cdot
\rho_{n+1}^{\indalp{i}{0}} A_1^{\indalp{i}{1}} \cdots A_{n-1}^{\indalp{i}{n-1}} \rho_i \cdot
\rho_n^{-1} A_{n-1}^{-\indalp{n}{n-1}} \cdots A_1^{-\indalp{n}{1}} \rho_{n+1}^{-\indalp{n}{0}}\\
=&
\rho_{n+1}^{\indalp{n}{0}- \indalp{i}{0}} A_1^{\epsilon(\indalp{i}{0}) \indalp{n}{1}} \cdots A_{n-1}^{\epsilon(\indalp{i}{0}) \indalp{n}{n-1}}
A_1^{-\delta(\indalp{i}{0})} \cdots A_{n-1}^{-\delta(\indalp{i}{0})} \rho_n \rho_{n+1}^{\indalp{n}{0}}\cdot\\
& A_1^{\epsilon(\indalp{n}{0}) \indalp{i}{1}} \cdots A_{n-1}^{\epsilon(\indalp{n}{0}) \indalp{i}{n-1}} 
A_i^{\delta(\indalp{n}{0})} \rho_i A_1^{\delta(\indalp{n}{0})}\cdots A_{n-1}^{\delta(\indalp{n}{0})} \rho_n^{-1} \cdot\\
& A_{n-1}^{-\epsilon(\indalp{n}{0}) \indalp{n}{n-1}} \cdots A_1^{-\epsilon(\indalp{n}{0}) \indalp{n}{1}}\\
=&
\rho_{n+1}^{-\indalp{i}{0}} A_1^{\epsilon(\indalp{n}{0}) \epsilon(\indalp{i}{0}) \indalp{n}{1}} \cdots A_{n-1}^{\epsilon(\indalp{n}{0}) \epsilon(\indalp{i}{0}) \indalp{n}{n-1}}
A_1^{-\epsilon(\indalp{n}{0}) \delta(\indalp{i}{0})} \cdots A_{n-1}^{-\epsilon(\indalp{n}{0}) \delta(\indalp{i}{0})}\cdot\\
& A_1^{-\delta(\indalp{n}{0})} \cdots A_{n-1}^{-\delta(\indalp{n}{0})} 
A_1^{\epsilon(\indalp{n}{0}) \indalp{i}{1}} \cdots A_{n-1}^{\epsilon(\indalp{n}{0}) \indalp{i}{n-1}} 
A_i^{\delta(\indalp{n}{0})} A_1^{\delta(\indalp{n}{0})}\cdots A_{n-1}^{\delta(\indalp{n}{0})} \cdot\\
& A_{n-1}^{-\epsilon(\indalp{n}{0}) \indalp{n}{n-1}} \cdots A_1^{-\epsilon(\indalp{n}{0}) \indalp{n}{1}} \rho_n \rho_i \rho_n^{-1}.
\end{align*}

Comparing coefficients in $K/L$, we obtain:
\begin{align}
& \text{$\indbet{i}{n}{0}=2\indalp{i}{0}$, so $\indbet{i}{n}{0}$ is even, $\epsilon(\indbet{i}{n}{0})=1$ and $\delta(\indbet{i}{n}{0})=0$}\label{eq:rhonrho1}\\
& \text{$\epsilon(\indalp{n}{0}) \indalp{n}{k}(\epsilon(\indalp{i}{0})-1)
+ \indalp{i}{k} (\epsilon(\indalp{n}{0}) -1) -\epsilon(\indalp{n}{0}) \delta(\indalp{i}{0}) = -\indbet{i}{n}{k}$ 
for $k=1,\ldots,n-1$, $k\neq i$}\notag\\
& \text{$\epsilon(\indalp{n}{0}) \indalp{n}{i}(\epsilon(\indalp{i}{0})-1)
+ \indalp{i}{i} (\epsilon(\indalp{n}{0}) -1) -\epsilon(\indalp{n}{0}) \delta(\indalp{i}{0}) +\delta(\indalp{n}{0}) = -\indbet{i}{n}{i}$.}\notag
\end{align}
In particular, the coefficient $\indbet{i}{n}{0}$ of $\rho_{n+1}$ in $\overline{s}_{\ast}(B_{i,n})$ is even. Further: 
\begin{align}
& \text{$\indbet{i}{n}{k}\equiv \delta(\indalp{i}{0}) \pmod 2$ for all $k\neq i$}\notag\\
& \text{$\indbet{i}{n}{i} \equiv \delta(\indalp{i}{0}) +\delta(\indalp{n}{0}) \pmod 2$ for all $1\leq i\leq n-1$.}\label{eq:ini}
\end{align}

\item Consider the relation \underline{$\rho_i^2= B_{1,i} \cdots B_{i-1,i} B_{i,i+1}\cdots B_{i,n}$} in $P_n(\rp)$, where $1\leq i\leq n-1$. Using equations~(\ref{eq:rho2rho1}) and~(\ref{eq:rhonrho1}), we see that:
\begin{align*}
\overline{s}_{\ast}(B_{1,i} \cdots B_{i-1,i} B_{i,i+1}\cdots B_{i,n-1} B_{i,n}) =& 
A_1^{\indbet{1}{i}{1}} \cdots A_{n-1}^{\indbet{1}{i}{n-1}} B_{1,i} \cdots
A_1^{\indbet{i-1}{i}{1}} \cdots A_{n-1}^{\indbet{i-1}{i}{n-1}} B_{i-1,i} \cdot\\
& A_1^{\indbet{i}{i+1}{1}} \cdots A_{n-1}^{\indbet{i}{i+1}{n-1}} B_{i,i+1} \cdots
A_1^{\indbet{i}{n-1}{1}} \cdots A_{n-1}^{\indbet{i}{n-1}{n-1}} \cdot\\
& B_{i,n-1} \rho_{n+1}^{2\indalp{i}{0}} A_1^{\indbet{i}{n}{1}} \cdots A_{n-1}^{\indbet{i}{n}{n-1}} B_{i,n}\\
=& 
\rho_{n+1}^{2\indalp{i}{0}} A_1^{\indbet{1}{i}{1}+\cdots +\indbet{i-1}{i}{1}+\indbet{i}{i+1}{1}+ \cdots+ \indbet{i}{n-1}{1}+ \indbet{i}{n}{1}}\cdot\\
& \cdots A_{n-1}^{\indbet{1}{i}{n-1}+ \cdots + \indbet{i-1}{i}{n-1}+\indbet{i}{i+1}{n-1}+\cdots + \indbet{i}{n-1}{n-1}+ \indbet{i}{n}{n-1}}\cdot\\
& B_{1,i} \cdots B_{i-1,i} B_{i,i+1} \cdots B_{i,n-1} B_{i,n}\\
=& 
\rho_{n+1}^{2\indalp{i}{0}} A_1^{\indbet{1}{i}{1}+\cdots +\indbet{i-1}{i}{1}+\indbet{i}{i+1}{1}+ \cdots+ \indbet{i}{n-1}{1}+ \indbet{i}{n}{1}}\cdot\\
& \cdots A_i^{\indbet{1}{i}{i}+ \cdots + \indbet{i-1}{i}{i}+\indbet{i}{i+1}{i}+\cdots + \indbet{i}{n-1}{i}+ \indbet{i}{n}{i}-1}\cdot\\
& \cdots A_{n-1}^{\indbet{1}{i}{n-1}+ \cdots + \indbet{i-1}{i}{n-1}+\indbet{i}{i+1}{n-1}+\cdots + \indbet{i}{n-1}{n-1}+ \indbet{i}{n}{n-1}} \rho_i^2,
\end{align*}
using the relation $B_{1,i} \cdots B_{i-1,i} B_{i,i+1} \cdots B_{i,n-1} B_{i,n} A_i=\rho_i^2$ in $P_{n+1}(\rp)/L$.

\begin{align*}
\overline{s}_{\ast}(\rho_i^2) =& \rho_{n+1}^{\indalp{i}{0}} A_1^{\indalp{i}{1}} \cdots A_{n-1}^{\indalp{i}{n-1}} \rho_i \cdot \rho_{n+1}^{\indalp{i}{0}} A_1^{\indalp{i}{1}} \cdots A_{n-1}^{\indalp{i}{n-1}} \rho_i\\
=& 
\rho_{n+1}^{2\indalp{i}{0}} A_1^{\indalp{i}{1}(\epsilon(\indalp{i}{0})+1)} \cdots A_i^{\indalp{i}{i}(\epsilon(\indalp{i}{0})+1)+ \delta(\indalp{i}{0})}
A_{n-1}^{\indalp{i}{n-1}(\epsilon(\indalp{i}{0})+1)} \rho_i^2.
\end{align*}

Comparing coefficients in $K/L$, for all $1\leq i\leq n-1$, we obtain:
\begin{align}
& \text{$\indbet{1}{i}{k}+ \cdots + \indbet{i-1}{i}{k}+\indbet{i}{i+1}{k}+\cdots + \indbet{i}{n-1}{k}+ \indbet{i}{n}{k}= \indalp{i}{k}(\epsilon(\indalp{i}{0})+1)$ for all $k\neq i$}\notag \\
& \indbet{1}{i}{i}+ \cdots + \indbet{i-1}{i}{i}+\indbet{i}{i+1}{i}+\cdots + \indbet{i}{n-1}{i}+ \indbet{i}{n}{i}-1=\indalp{i}{i}(\epsilon(\indalp{i}{0})+1)+ \delta(\indalp{i}{0}).\label{eq:rhoisq}
\end{align}

\item Consider the relation \underline{$\rho_n^2= B_{1,n} \cdots B_{n-1,n}$}  in $P_n(\rp)$:
\begin{align*}
\overline{s}_{\ast}(B_{1,n} \cdots B_{n-1,n}) =& 
\rho_n^{2\indalp{1}{0}} A_1^{\indbet{1}{n}{1}} \cdots A_{n-1}^{\indbet{1}{n}{n-1}} B_{1,n} \cdots 
\rho_n^{2\indalp{n-1}{0}} A_1^{\indbet{n-1}{n}{1}} \cdots A_{n-1}^{\indbet{n-1}{n}{n-1}} B_{n-1,n}\\
=& \rho_n^{2(\indalp{1}{0}+\cdots + \indalp{n-1}{0})} A_1^{\indbet{1}{n}{1}+\cdots + \indbet{n-1}{n}{1}} \cdots A_{n-1}^{\indbet{1}{n}{n-1}+\cdots + \indbet{n-1}{n}{n-1}}\cdot\\
& B_{1,n} \cdots B_{n-1,n}\\
=& \rho_n^{2(\indalp{1}{0}+\cdots + \indalp{n-1}{0})} A_1^{\indbet{1}{n}{1}+\cdots + \indbet{n-1}{n}{1}} \cdots A_{n-1}^{\indbet{1}{n}{n-1}+\cdots + \indbet{n-1}{n}{n-1}}\cdot\\
& \rho_n^2 \rho_{n+1}^{-2} A_1 \cdots A_{n-1}\\
=& \rho_n^{2(\indalp{1}{0}+\cdots + \indalp{n-1}{0}-1)} A_1^{\indbet{1}{n}{1}+\cdots + \indbet{n-1}{n}{1}+1} \cdots A_{n-1}^{\indbet{1}{n}{n-1}+\cdots + \indbet{n-1}{n}{n-1}+1} \rho_n^2, \end{align*}
using the relations $B_{1,n} \cdots B_{n-1,n} B_{n,n+1}=\rho_n^2$ and $A_1 \cdots A_{n-1} B_{n,n+1}=\rho_{n+1}^2$, and the fact that $\rho_n^2 \comm \rho_{n+1}^2$ in $P_{n+1}(\rp)/L$.

\begin{align*}
\overline{s}_{\ast}(\rho_n^2) =&
\rho_{n+1}^{\indalp{n}{0}} A_1^{\indalp{n}{1}} \cdots A_{n-1}^{\indalp{n}{n-1}} \rho_n \cdot \rho_{n+1}^{\indalp{n}{0}} A_1^{\indalp{n}{1}} \cdots A_{n-1}^{\indalp{n}{n-1}} \rho_n\\
=&
\rho_{n+1}^{\indalp{n}{0}} A_1^{\indalp{n}{1}} \cdots A_{n-1}^{\indalp{n}{n-1}} 
\rho_{n+1}^{-\indalp{n}{0}} A_1^{-\delta(\indalp{n}{0})} \cdots A_{n-1}^{-\delta(\indalp{n}{0})} A_1^{\indalp{n}{1}} \cdots A_{n-1}^{\indalp{n}{n-1}} \rho_n^2\\
=&
A_1^{\epsilon(\indalp{n}{0}) \indalp{n}{1}} \cdots A_{n-1}^{\epsilon(\indalp{n}{0}) \indalp{n}{n-1}} 
A_1^{-\delta(\indalp{n}{0})} \cdots A_{n-1}^{-\delta(\indalp{n}{0})} A_1^{\indalp{n}{1}} \cdots A_{n-1}^{\indalp{n}{n-1}} \rho_n^2\\
=&
A_1^{\indalp{n}{1}(\epsilon(\indalp{n}{0})+1) -\delta(\indalp{n}{0})} \cdots
A_{n-1}^{\indalp{n}{n-1}(\epsilon(\indalp{n}{0})+1) -\delta(\indalp{n}{0})}
\rho_n^2.
\end{align*}

Comparing coefficients in $K/L$, we obtain:
\begin{align}
& \indalp{1}{0}+\cdots + \indalp{n-1}{0}=1\label{eq:sumone}\\
& \text{$\indbet{1}{n}{i}+\cdots + \indbet{n-1}{n}{i}+1= \indalp{n}{i}(\epsilon(\indalp{n}{0})+1) -\delta(\indalp{n}{0})$ for all $i=1,\ldots, n-1$.}\notag
\end{align}
\end{enumerate}

Now consider equation~(\ref{eq:rhoisq}) modulo~$2$. For all $1\leq i\leq n-1$, we have:
\begin{align*}
\delta(\indalp{i}{0}) & \equiv \indbet{1}{i}{i}+ \cdots + \indbet{i-1}{i}{i}+\indbet{i}{i+1}{i}+\cdots + \indbet{i}{n-1}{i}+ \indbet{i}{n}{i}+1\\
& \equiv \delta(\indalp{1}{0})+ \cdots +\delta(\indalp{i-1}{0}) + \delta(\indalp{i+1}{0}) + \cdots +\delta(\indalp{n-1}{0})+ (\delta(\indalp{i}{0})+\delta(\indalp{n}{0}))+1,
\end{align*}
using equations~(\ref{eq:iji}),~(\ref{eq:ijj}) and~(\ref{eq:ini}). Hence 
\begin{equation}\label{eq:deltasum}
\text{$\delta(\indalp{i}{0}) \equiv 1+\sum_{j=1}^n \; \delta(\indalp{j}{0}) \pmod 2$ for all $1\leq i\leq n-1$,}
\end{equation}
and thus $\delta(\indalp{1}{0}) \equiv \cdots \equiv \delta(\indalp{n-1}{0}) \pmod 2$. Further, since $x \equiv \delta(x) \pmod 2$ for all $x\in \Z$, we see from equation~(\ref{eq:sumone}) that $\sum_{j=1}^{n-1} \; \delta(\indalp{j}{0}) \equiv 1 \pmod 2$, $\delta(\indalp{1}{0}) \equiv \cdots \equiv \delta(\indalp{n-1}{0})\equiv 1 \pmod 2$ and that $n$ is even. It follows from equation~(\ref{eq:deltasum}) that $\delta(\indalp{1}{0}) \equiv \cdots \equiv \delta(\indalp{n-1}{0}) \equiv \delta(\indalp{n}{0}) \equiv 1\pmod 2$, and so $\indalp{1}{0}, \ldots, \indalp{n}{0}$ are odd. Since $n$ is even, the element $B_{2,3}$ exists. Further, $3\leq n-1$, and hence $\indbet{2}{3}{0}=0$ from equation~(\ref{eq:rho2rho1}). Now consider the image in $P_{n+1}(\rp)/L$ under $\overline{s}_{\ast}$ of the relation $\rho_1 \comm B_{2,3}$ of $P_n(\rp)$: 
\begin{align*}
\overline{s}_{\ast}(\rho_1 B_{2,3}) &= \rho_{n+1}^{\indalp{1}{0}} A_1^{\indalp{1}{1}} \cdots A_{n-1}^{\indalp{1}{n-1}} \rho_1 \cdot  A_1^{\indbet{2}{3}{1}} \cdots A_{n-1}^{\indbet{2}{3}{n-1}} B_{2,3}\\
&= \rho_{n+1}^{\indalp{1}{0}} A_1^{\indalp{1}{1}+ \indbet{2}{3}{1}} \cdots A_{n-1}^{\indalp{1}{n-1}+ \indbet{2}{3}{n-1}} \rho_1 B_{2,3}.\\
\overline{s}_{\ast}(B_{2,3} \rho_1) &= \rho_{n+1}^{\indalp{1}{0}} A_1^{\epsilon(\indalp{1}{0}) \indbet{2}{3}{1}+ \indalp{1}{1}} \cdots A_{n-1}^{\epsilon(\indalp{1}{0}) \indbet{2}{3}{n-1}+ \indalp{1}{n-1}} B_{2,3} \rho_1.
\end{align*}

Comparing coefficients in $K/L$, we see that $\indbet{2}{3}{i}(\epsilon(\indalp{1}{0})-1)=0$ for all $i=1,\ldots, n-1$. Since $\indalp{1}{0}$ is odd, $\epsilon(\indalp{1}{0})=-1$, and thus $\indbet{2}{3}{i}=0$ for all $i=1,\ldots, n-1$. Hence $\indbet{2}{3}{2}=0$. But since $n$ is even, $3\leq n-1$, and this contradicts equation~(\ref{eq:iji}). Hence $\overline{p}_{\ast}$ does not admit a section, and so neither does $p_{\ast}$. This proves the first statement of the theorem. The second statement follows from the fact that we mentioned in the introduction, that under the hypotheses of the theorem, the fibration $\map{p}{F_{n+m}(\rp)}[F_n(\rp)]$ admits a section if and only if the group homomorphism $\map{p_{\ast}}{P_{n+m}(\rp)}[P_n(\rp)]$ does.
\end{proof}

\end{document}